\documentclass[10pt,twoside]{article}
\usepackage{graphicx}
\usepackage{amsmath}
\usepackage{Latex-document}
\def\Section#1{\section{\hskip -1em . \hskip 0.8em #1}}

\title{\bf  Permutation Groups and \vskip -2mm Normal Subgroups\vskip 6mm}
\author{Cheryl E. Praeger\vspace*{-0.5cm}\thanks{Department of Mathematics
\& Statistics, University of Western Australia, 35 Stirling
Highway, Crawley, Western Australia 6009, Australia. E-mail:
praeger@maths.uwa.edu.au}}
\date{\vspace{-8mm}}

\begin{document}

\maketitle

\newcommand{\Sym}{{\rm Sym}}
\newcommand{\Aut}{{\rm Aut}}
\newcommand{\Comp}{{\sf Comp}}
\newcommand{\Om}{\Omega}
\newcommand{\De}{\Delta}
\newcommand{\al}{\alpha}
\newcommand{\be}{\beta}
\newcommand{\de}{\delta}
\newcommand{\ga}{\gamma}
\newcommand{\calP}{{\mathcal{P}}}
\newcommand{\soc}{{\rm soc}}
\renewcommand{\thesection}{\arabic{section}}
\def\addsec{\setcounter{theorem}{0}}
\newtheorem{theorem}{Theorem}[section]

\thispagestyle{first} \setcounter{page}{67}

\begin{abstract}

\vskip 3mm

Various descending chains of subgroups of a finite permutation
group can be used to define a sequence of `basic' permutation
groups that are analogues of composition factors for abstract
finite groups.  Primitive groups have been the traditional choice
for this purpose, but some combinatorial applications require
different kinds of basic groups, such as quasiprimitive groups,
that are defined by properties of their normal subgroups.
Quasiprimitive groups admit similar analyses to primitive groups,
share many of their properties, and have been used successfully,
for example to study $s$-arc transitive graphs.  Moreover
investigating them has led to new results about finite simple
groups.

\vskip 4.5mm


\end{abstract}

\vskip 12mm

\Section{Introduction} \label{section 1}\setzero\addsec

\vskip-5mm \hspace{5mm}

For a satisfactory understanding of finite groups it is important
to study simple groups and characteristically simple groups, and
how to fit them together to form arbitrary finite groups.  This
paper discusses an analogous programme for studying finite
permutation groups.  By considering various descending subgroup
chains of finite permutation groups we define in \S2 sequences of
`basic' permutation groups that play the role for finite
permutation groups that composition factors or chief factors play
for abstract finite groups.  Primitive groups have been the
traditional choice for basic permutation groups, but for some
combinatorial applications larger families of basic groups, such
as quasiprimitive groups, are needed (see \S3).

Application of a theorem first stated independently in 1979 by M.
E. O'Nan and L. L. Scott~\cite{Cameron81} has proved to be the
most useful modern method for identifying the possible structures
of finite primitive groups, and is now used routinely for their
analysis.  Analogues of this theorem are available for the
alternative families of basic permutation groups. These theorems
have become standard tools for studying finite combinatorial
structures such as vertex-transitive graphs and examples are given
in \S3 of successful analyses for distance transitive graphs and
$s$-arc-transitive graphs. Some characteristic properties of basic
permutation groups, including these structure theorems are
discussed in \S4.

Studying the symmetry of a family of finite algebraic or
combinatorial systems often leads to problems about groups of
automorphisms acting as basic permutation groups on points or
vertices. In particular determining the full automorphism group of
such a system sometimes requires a knowledge of the permutation
groups containing a given basic permutation group, and for this it
is important to understand the lattice of basic permutation groups
on a given set. The fundamental problem here is that of
classifying all inclusions of one basic permutation group in
another, and integral to its solution is a proper understanding of
the factorisations of simple and characteristically simple groups.
In \S3 and \S4 we outline the current status of our knowledge
about such inclusions and their use.

The precision of our current knowledge of basic permutation groups
depends heavily on the classification of the finite simple groups.
Some problems about basic permutation groups translate directly to
questions about simple groups, and answering them leads to new
results about simple groups. Several of these results and their
connections with basic groups are discussed in the final section
\S5.

In summary, this approach to analysing finite permutation groups
involves an interplay between combinatorics, group actions, and
the theory of finite simple groups.  One measure of its success is
its effectiveness in combinatorial applications.

\Section{Defining basic permutation groups} \label{section 2} \setzero
\addsec

\vskip-5mm \hspace{5mm}

Let $G$ be a subgroup of the symmetric group $\Sym(\Om)$ of all
permutations of a finite set $\Om$.  Since an intransitive
permutation group is contained in the direct product of its
transitive constituents, it is natural when studying permutation
groups to focus first on the transitive ones.  Thus we will assume
that $G$ is transitive on $\Om$.  Choose a point $\al\in\Om$ and
let $G_\al$ denote the subgroup of $G$ of permutations that fix
$\al$, that is, the stabiliser of $\al$.  Let ${\sf Sub}(G,G_\al)$
denote the lattice of subgroups of $G$ containing $G_\al$.  The
concepts introduced below are independent of the choice of $\al$
because of the transitivity of $G$.  We shall introduce three
types of basic permutation groups, relative to
$\mathcal{L}_1:={\sf Sub}(G,G_\al)$ and two other types of
lattices $\mathcal{L}_2$ and $\mathcal{L}_3$, where we regard each
$\mathcal{L}_i$ as a function that can be evaluated on any finite
transitive group $G$ and stabiliser $G_\al$.

For $G_\al\leq H\leq G$, the $H$-orbit containing $\al$ is
$\al^H=\{\al^h\,|\,h\in H\}$.  If $G_\al\leq H<K\leq G$, then the
$K$-images of $\al^H$ form the parts of a $K$-invariant partition
$\calP(K,H)$ of $\al^K$, and $K$ induces a transitive permutation
group $\Comp(K,H)$ on $\calP(K,H)$ called a \emph{component} of
$G$.  In particular the component $\Comp(G,G_\al)$ permutes
$\calP(G,G_\al)=\{\{\beta\}\,|\,\beta\in\Om\}$ in the same way
that $G$ permutes $\Om$, and we may identify $G$ with
$\Comp(G,G_\al)$.

For a lattice $\mathcal{L}$ of subgroups of $G$ containing
$G_\alpha$, we say that $K$ \emph{covers $H$ in $\mathcal{L}$} if
$K, H \in \mathcal{L}$, $H<K$, and there are no intermediate
subgroups lying in $\mathcal{L}$.  The \emph{basic components of
$G$ relative to $\mathcal{L}$} are then defined as all the
components $\Comp(K,H)$ for which $K$ covers $H$ in $\mathcal{L}$.
Each maximal chain $G_\al=G_0<G_1<\dots<G_r=G$ in $\mathcal{L}$
determines a sequence of basic components relative to
$\mathcal{L}$, namely
$\Comp(G_1,G_{0}),\dots,\Comp(G_{r},G_{r-1})$, and $G$ can be
embedded in the iterated wreath product of these groups.  In this
way the permutation groups occurring as basic components relative
to $\mathcal{L}$, for some finite transitive group, may be
considered as `building blocks' for finite permutation groups.  We
refer to such groups as basic permutation groups relative to
$\mathcal{L}$.

A transitive permutation group $G$ on $\Om$ is \emph{primitive} if
$G_\al$ is a maximal subgroup of $G$, that is, if ${\sf
Sub}(G,G_\al)=\{G,G_\al\}$.  The basic components of $G$ relative
to $\mathcal{L}_1={\sf Sub}(G,G_\al)$ are precisely those of its
components that are primitive.

The basic groups of the second type are the quasiprimitive groups.
A transitive permutation group $G$ on $\Om$ is
\emph{quasiprimitive} if each nontrivial normal subgroup of $G$ is
transitive on $\Om$.  The corresponding sublattice is the set
$\mathcal{L}_2$ of all subgroups $H\in {\sf Sub}(G,G_\al)$ such
that there is a sequence $H_0=H\leq H_1\leq \dots\leq H_r=G$ with
each subgroup of the form $H_i=G_\al N_i$ where for $i<r$, $N_i$
is a normal subgroup of $H_{i+1}$, and $N_r=G$.  The basic
components of $G$ relative to $\mathcal{L}_2$ are precisely those
of its components that are quasiprimitive.

Basic groups of the third type are \emph{innately transitive},
namely transitive permutation groups that have at least one
transitive minimal normal subgroup. The corresponding sublattice
will be $\mathcal{L}_3$.  A subgroup $N$ of $G$ is
\emph{subnormal} in $G$ if there is a sequence $N_0=N\leq N_1\leq
\dots\leq N_r=G$ such that, for $i<r$, $N_i$ is a normal subgroup
of $N_{i+1}$.  The lattice $\mathcal{L}_3$ consists of all
subgroups of the form $G_\al N$, where $N$ is subnormal in $G$ and
normalised by $G_\al$.  All the basic components of $G$ relative
to $\mathcal{L}_3$ are innately transitive.  Note that each
primitive group is quasiprimitive and each quasiprimitive group is
innately transitive.  Proofs of the assertions about
$\mathcal{L}_2$ and $\mathcal{L}_3$ and their components may be
found in \cite{seminorm}.

\Section{The role of basic groups in graph theory} \setzero \addsec

\vskip-5mm \hspace{5mm}

For many group theoretic and combinatorial applications finite
primitive permutation groups are the appropriate basic permutation
groups, since many problems concerning finite permutation groups
can be reduced to the case of primitive groups. However such
reductions are sometimes not possible when studying
point-transitive automorphism groups of finite combinatorial
structures because the components of the given point-transitive
group have no interpretation as point-transitive automorphism
groups of structures within the family under investigation.  The
principal motivation for studying some of these alternative basic
groups came from graph theory, notably the study of $s$-arc
transitive graphs ($s\geq2$).

A finite graph $\Gamma=(\Om,E)$ consists of a finite set $\Om$ of
points, called vertices, and a subset $E$ of unordered pairs from
$\Om$ called edges.  For $s\geq1$, an \emph{$s$-arc} of $\Gamma$
is a vertex sequence $(\al_0,\al_1,\dots,\al_s)$ such that each
$\{\al_i,\al_{i+1}\}$ is an edge and $\al_{i-1}\ne\al_{i+1}$ for
all $i$.  We usually call a 1-arc simply an arc. Automorphisms of
$\Gamma$ are permutations of $\Om$ that leave $E$ invariant, and a
subgroup $G$ of the automorphism group ${\rm Aut}(\Gamma)$ is
\emph{$s$-arc-transitive} if $G$ is transitive on the $s$-arcs of
$\Gamma$.  If $\Gamma$ is connected and is regular of
\emph{valency} $k>0$ so that each vertex is in $k$ edges, then an
$s$-arc-transitive subgroup $G\leq{\rm Aut}(\Gamma)$ is in
particular transitive on $\Om$ and also, if $s\geq2$, on
$(s-1)$-arcs.  It is natural to ask which of the components of
this transitive permutation group $G$ on $\Om$ act as
$s$-arc-transitive automorphism groups of graphs related to
$\Gamma$.

For $G_\al\leq H\leq G$, there is a naturally defined quotient
graph $\Gamma_H$ with vertex set the partition of $\Om$ formed by
the $G$-images of the set $\al^H$, where two such $G$-images are
adjacent in $\Gamma_H$ if at least one vertex in the first is
adjacent to at least one vertex of the second.  If $\Gamma$ is
connected and $G$ is arc-transitive, then $\Gamma_H$ is connected
and $G$ induces an arc-transitive automorphism group of
$\Gamma_H$, namely the component $\Comp(G,H)$.  If $H$ is a
maximal subgroup of $G$, then $\Comp(G,H)$ is both
vertex-primitive and arc-transitive on $\Gamma_H$.  This
observation enables many questions about arc-transitive graphs to
be reduced to the vertex-primitive case.

Perhaps the most striking example is provided by the family of
finite distance transitive graphs. The \emph{distance} between two
vertices is the minimum number of edges in a path joining them,
and $G$ is \emph{distance transive} on $\Gamma$ if for each $i$,
$G$ is transitive on the set of ordered pairs of vertices at
distance $i$.  In particular if $G$ is distance transitive on
$\Gamma$ then $\Gamma$ is connected and regular, of valency $k$
say.  If $k=2$ then $\Gamma$ is a cycle and all cycles are
distance transitive, so suppose that $k\geq3$.  If $\Gamma_H$ has
more than two vertices, then $\Comp(G,H)$ is distance transitive
on $\Gamma_H$, while if $\Gamma_H$ has only two vertices then $H$
is distance transitive on a smaller graph $\Gamma_2$, namely
$\Gamma_2$ has $\al^H$ as vertex set with two vertices adjacent if
and only if they are at distance 2 in $\Gamma$ (see for example
\cite{sasha}).  Passing to $\Gamma_H$ or $\Gamma_2$ respectively
and repeating this process, we reduce to a vertex-primitive
distance transitive graph.  The programme of classifying the
finite vertex-primitive distance transitive graphs is approaching
completion, and surveys of progress up to the mid 1990's can be
found in \cite{sasha,vbc}.  The initial result that suggested a
classification might be possible is the following.  Here a group
$G$ is \emph{almost simple} if $T\leq G\leq {\rm Aut}(T)$ for some
nonabelian simple group $T$, and a permutation group $G$ has
\emph{affine type} if $G$ has an elementary abelian regular normal
subgroup.

\begin{theorem}\label{dtg}{\rm\cite{psy}} If $G$ is
vertex-primitive and distance transitive on a finite graph
$\Gamma$, then either $\Gamma$ is known explicitly, or $G$ is
almost simple, or $G$ has affine type.
\end{theorem}

In general, if $G$ is $s$-arc-transitive on $\Gamma$ with
$s\geq2$, then none of the components $\Comp(G,H)$ with
$G_\al<H<G$ is $s$-arc-transitive on $\Gamma_H$, so there is no
hope that the problem of classifying finite $s$-arc-transitive
graphs, or even giving a useful description of their structure,
can be reduced to the case of vertex-primitive $s$-arc-transitive
graphs.  However the class of $s$-arc transitive graphs behaves
nicely with respect to \emph{normal quotients}, that is, quotients
$\Gamma_H$ where $H=G_\al N$ for some normal subgroup $N$ of $G$.
For such quotients, the vertex set of $\Gamma_H$ is the set of
$N$-orbits, $G$ acts $s$-arc-transit\-ive\-ly on $\Gamma_H$, and
if $\Gamma_H$ has more than two vertices then $\Gamma$ is a cover
of $\Gamma_H$ in the sense that, for two $N$-orbits adjacent in
$\Gamma_H$, each vertex in one $N$-orbit is adjacent in $\Gamma$
to exactly one vertex in the other $N$-orbit. We say that $\Gamma$
is a \emph{normal cover} of $\Gamma_H$. If in addition $N$ is a
maximal intransitive normal subgroup of $G$ with more than two
orbits, then $G$ is both vertex-quasiprimitive and
$s$-arc-transitive on $\Gamma_H$, see \cite{2arcqp}.  If some
quotient $\Gamma_H$ has two vertices then $\Gamma$ is bipartite,
and such graphs require a specialised analysis that parallels the
one described here. On the other hand if $\Gamma$ is not bipartite
then $\Gamma$ is a normal cover of at least one $\Gamma_H$ on
which the $G$-action is both vertex-quasiprimitive and
$s$-arc-transitive.  The wish to understand quasiprimitive $s$-arc
transitive graphs led to the development of a theory for finite
quasiprimitive permutation groups similar to the theory of finite
primitive groups.  Applying this theory led to a result similar to
Theorem~\ref{dtg}, featuring two additional types of
quasiprimitive groups, called \emph{twisted wreath type} and
\emph{product action type}. Descriptions of these types may be
found in \cite{2arcqp} and \cite{bcc}.

\begin{theorem}\label{2arc}{\rm \cite{2arcqp}}
If $G$ is vertex-quasiprimitive and $s$-arc-transitive on a finite
graph $\Gamma$ with $s\geq2$, then $G$ is almost simple, or of
affine, twisted wreath or product action type.
\end{theorem}

Examples exist for each of the four quasiprimitive types, and
moreover this division of vertex-quasiprimitive $s$-arc transitive
graphs into four types has resulted in a better understanding of
these graphs, and in some cases complete classifications.  For
example all examples with $G$ of affine type, or with $T\leq
G\leq{\Aut}(T)$ and $T={\rm PSL}_2(q), {\rm Sz}(q)$ or ${\rm
Ree}(q)$ have been classified, in each case yielding new $s$-arc
transitive graphs, see \cite{ip,bcc}.  Also using
Theorem~\ref{2arc} to study the normal quotients of an $s$-arc
transitive graph has led to some interesting restrictions on the
number of vertices.

\begin{theorem}\label{li}{\rm \cite{li1,li2}}
Suppose that $\Gamma$ is a finite $s$-arc-transitive graph with
$s\geq4$. Then the number of vertices is even and not a power of
$2$.
\end{theorem}

The concept of a normal quotient has proved useful for analysing
many families of edge-transitive graphs, even those for which a
given edge-transitive group is not vertex-transitive.  For example
it provides a framework for a systematic study of locally
$s$-arc-transitive graphs in which quasiprimitive actions are of
central importance, see \cite{giudici2}.

We have described how to form primitive arc-transitive quotients
of arc-trans\-it\-ive graphs, and quasiprimitive
$s$-arc-transitive normal quotients of non-bipartite
$s$-arc-transitive graphs.  However recognising these quotients is
not always easy without knowing their full automorphism groups.
To identify the automorphism group of a graph, given a primitive
or quasiprimitive subgroup $G$ of automorphisms, it is important
to know the permutation groups of the vertex set that contain $G$,
that is the over-groups of $G$.  In the case of finite primitive
arc-transitive and edge-transitive graphs, knowledge of the
lattice of primitive permutation groups on the vertex set together
with detailed knowledge of finite simple groups led to the
following result.  The \emph{socle} of a finite group $G$, denoted
$\soc(G)$, is the product of its minimal normal subgroups.

\begin{theorem}\label{aut}{\rm \cite{LiebeckPraegeretal02}}
Let $G$ be a primitive arc- or edge-transitive group of
automorphisms of a finite connected graph $\Gamma$.  Then either
$G$ and ${\rm Aut}(\Gamma)$ have the same socle, or $G < H\leq
{\rm Aut}(\Gamma)$ where $\soc(G)\ne \soc(H)$ and $G, H$ are
explicitly listed.
\end{theorem}

In the case of graphs $\Gamma$ for which a quasiprimitive subgroup
$G$ of $\Aut(\Gamma)$ is given, it is possible that $\Aut(\Gamma)$
may not be quasiprimitive.  However, even in this case a good
knowledge of the quasiprimitive over-groups of a quasiprimitive
group is helpful, for if $N$ is a maximal intransitive normal
subgroup of $\Aut(\Gamma)$ then both $G$ and $\Aut(\Gamma)$ induce
quasiprimitive automorphism groups of the normal quotient
$\Gamma_H$, where $H=\Aut(\Gamma)_\alpha N$, and the action of $G$
is faithful. This approach was used, for example, in classifying
the $2$-arc transitive graphs admitting ${\rm Sz}(q)$ or ${\rm
Ree}(q)$ mentioned above, and also in analysing the automorphism
groups of Cayley graphs of simple groups in \cite{cayley}.

Innately transitive groups, identified in \S\,\ref{section 2} as a
third possibility for basic groups, have not received much
attention until recently. They arise naturally when investigating
the full automorphism groups of graphs. One example is given in
\cite{fhp} for locally-primitive graphs $\Gamma$ admitting an
almost simple vertex-quasiprimitive subgroup $G$ of automorphisms.
It is shown that either ${\rm Aut}(\Gamma)$ is innately
transitive, or $G$ is of Lie type in characteristic $p$ and ${\rm
Aut}(\Gamma)$ has a minimal normal $p$-subgroup involving a known
$G$-module.

\Section{Characteristics of basic permutation groups} \setzero \addsec

\vskip-5mm \hspace{5mm}

Finite primitive permutation groups have attracted the attention
of mathematicians for more than a hundred years.  In particular,
one of the central problems of 19th century Group Theory was to
find an upper bound, much smaller than $n!$, for the order of a
primitive group on a set of size $n$, other than the symmetric
group $S_n$ and the alternating group $A_n$.  It is now known that
the largest such groups occur for $n$ of the form $c(c-1)/2$ and
are $S_c$ and $A_c$ acting on the unordered pairs from a set of
size $c$.  The proofs of this and other results in this section
depend on the finite simple group classification.

If $G$ is a quasiprimitive permutation group on $\Om$,
$\al\in\Om$, and $H$ is a maximal subgroup of $G$ containing
$G_\al$, then the primitive component $\Comp(G,H)$ is isomorphic
to $G$ since the kernel of this action is an intransitive normal
subgroup of $G$ and hence is trivial.  Because of this we may
often deduce information about quasiprimitive groups from their
primitive components, and indeed it was found in
\cite{PraegerShalev01} that finite quasiprimitive groups possess
many characteristics similar to those of finite primitive groups.
This is true also of innately transitive groups.  We state just
one example, concerning the orders of permutation groups acting on
a set of size $n$, that is, of \emph{degree} $n$.

\begin{theorem}{\rm \cite{Cameron81,PraegerShalev01}}
There is a constant $c$ and an explicitly defined family
${\mathcal{F}}$ of finite permutation groups such that, if $G$ is
a primitive, quasiprimitive, or innately transitive permutation
group of degree $n$, then either $G\in{\mathcal{F}}$, or
$|G|<n^{c\log n}$.
\end{theorem}

The O'Nan-Scott Theorem partitions the finite primitive permutation groups into several disjoint types according
to the structure or action of their minimal normal subgroups.  It highlights the role of simple groups and their
representations in analysing and using primitive groups.  One of its first successful applications was the
analysis of distance transitive graphs in Theorem~\ref{dtg}.  Other early applications include a proof
\cite{CameronPraegeretal} of the Sims Conjecture, and a classification result~\cite{maxan} for maximal subgroups
of $A_n$ and $S_n$, both of which are stated below.

\begin{theorem}{\rm\cite{CameronPraegeretal}}
There is a function $f$ such that if $G$ is primitive on a finite
set $\Om$, and for $\al\in\Om$, $G_\al$ has an orbit of length $d$
in $\Om\setminus\{\al\}$, then $|G_\al|\leq f(d)$.
\end{theorem}

\begin{theorem}\label{maxan}{\rm \cite{maxan}}
Let $G=A_n$ or $S_n$ with $M$ a maximal subgroup.  Then either $M$
belongs to an explicit list or $M$ is almost simple and primitive.
Moreover if $H<G$ and $H$ is almost simple and primitive but not
maximal, then $(H,n)$ is known.
\end{theorem}

This is a rather curious way to state a classification result.
However it seems almost inconceivable that the finite almost
simple primitive groups will ever be listed explicitly.  Instead
\cite{maxan} gives an explicit list of triples $(H,M,n)$, where
$H$ is primitive of degree $n$ with a nonabelian simple normal
subgroup $T$ not normalised by $M$, and $H<M<HA_n$.  This result
suggested the possibility of describing the lattice of all
primitive permutation groups on a given set, for it gave a
description of the over-groups of the almost simple primitive
groups. Such a description was achieved in \cite{inclns} using a
general construction for primitive groups called a \emph{blow-up}
construction introduced by Kovacs \cite{kovacs}. The analysis
leading to Theorem~\ref{aut} was based on this theorem.

\begin{theorem}{\rm \cite{inclns}}
All inclusions $G<H<S_n$ with $G$ primitive are either explicitly
described, or are described in terms of a blow-up of an explicitly
listed inclusion $G_1<H_1<S_{n_1}$ with $n$ a proper power of
$n_1$.
\end{theorem}

Analogues of the O'Nan-Scott Theorem for finite quasiprimitive and innately transitive groups have been proved in
\cite{BambergPraeger,2arcqp} and enable similar analyses to be undertaken for problems involving these classes of
groups. For example, the quasiprimitive version formed the basis for Theorems~\ref{2arc} and \ref{li}.  It seems
to be the most useful version for dealing with families of vertex-transitive or locally-transitive graphs.  A
description of the lattice of quasiprimitive subgroups of $S_n$ was given in \cite{qpinclns,qpinclns2} and was
used, for example, in analysing Cayley graphs of finite simple groups in \cite{cayley}.

\begin{theorem}{\rm \cite{qpinclns,qpinclns2}}\label{qpinclns}
Suppose that $G<H<S_n$ with $G$ quasiprimitive and imprimitive, and $H$ quasiprimitive but $H\ne A_n$. Then either
$G$ and $H$ have equal socles and the same O'Nan-Scott types, or the possibilities for the O'Nan-Scott types of
$G, H$ are restricted and are known explicitly.
\end{theorem}

In the latter case, for most pairs of O'Nan-Scott types, explicit constructions are given for these inclusions.
Not all the types of primitive groups identified by the O'Nan-Scott Theorem occur for every degree $n$.  Let us
call permutation groups of degree $n$ other than $A_n$ and $S_n$ nontrivial.  A systematic study by Cameron,
Neumann and Teague \cite{CNT} of the integers $n$ for which there exists a nontrivial primitive group of degree
$n$ showed that the set of such integers has density zero in the natural numbers.  Recently it was shown in
\cite{PraegerShalev02} that a similar result holds for the degrees of nontrivial quasiprimitive and innately
transitive permutation groups. Note that $2.2< \sum_{d=1}^\infty\frac{1} {d\phi(d)}<2.23$.

\begin{theorem}{\rm \cite{CNT,PraegerShalev02}} \label{degs}
For a positive real number $x$, the proportion of integers $n\leq
x$ for which there exists a nontrivial primitive, quasiprimitive,
or innately transitive permutation group of degree $n$ is at most
$(1+o(1))c/\log x$, where $c=2$ in the case of primitive groups,
or $c=1+\sum_{d=1}^\infty\frac{1} {d\phi(d)}$ for the other cases.
\end{theorem}

\Section{Simple groups and basic permutation groups} \setzero \addsec

\vskip-5mm \hspace{5mm}

Many of the results about basic permutation groups mentioned above
rely on specific knowledge about finite simple groups.  Sometimes
this knowledge was already available in the simple group
literature.  However investigations of basic permutation groups
often raised interesting new questions about simple groups.
Answering these questions became an integral part of the study of
basic groups, and the answers enriched our understanding of finite
simple groups.  In this final section we review a few of these new
simple group results.  Handling the primitive almost simple
classical groups was the most difficult part of proving
Theorem~\ref{maxan}, and the following theorem of Aschbacher
formed the basis for their analysis.

\begin{theorem}\label{asch}{\rm\cite{asch}}
Let $G$ be a subgroup of a finite almost simple classical group
$X$ such that $G$ does not contain $\soc(X)$, and let $V$ denote
the natural vector space associated with $X$.  Then either $G$
lies in one of eight explicitly defined families of subgroups, or
$G$ is almost simple, absolutely irreducible on $V$ and the
(projective) representation of $\soc(G)$ on $V$ cannot be realised
over a proper subfield.
\end{theorem}

A detailed study of classical groups based on Theorem~\ref{asch}
led to Theorem~\ref{factns}, a classification  of the maximal
factorisations of the almost simple groups.  This classification
was fundamental to the proofs of Theorems~\ref{aut} and
\ref{maxan}, and has been used in diverse applications, for
example see \cite{fgs,li3}.

\begin{theorem}\label{factns}{\rm \cite{factns, max+}}
Let $G$ be a finite almost simple group and suppose that $G=AB$,
where $A, B$ are both maximal in $G$ subject to not containing
$\soc(G)$. Then $G, A, B$ are explicitly listed.
\end{theorem}

For a finite group $G$, let $\pi(G)$ denote the set of prime
divisors of $|G|$.  For many simple groups $G$ there are small
subsets of $\pi(G)$ that do not occur in the order of any proper
subgroup, and it is possible to describe some of these precisely
as follows.

\begin{theorem}\label{trans} {\rm \cite[Theorem 4, Corollaries 5 and 6]
{trans}} Let $G$ be an almost simple group with socle $T$, and let
$M$ be a subgroup of $G$ not containing $T$.
\begin{enumerate}
\item[(a)] If $G=T$ then for an explicitly defined subset $\Pi\subseteq
\pi(T)$ with $|\Pi|\leq3$, if $\Pi\subseteq\pi(M)$ then $T, M$ are
known explicitly, and in most cases $\pi(T)=\pi(M)$.
\item[(b)] If $\pi(T)\subseteq\pi(M)$ then $T, M$ are known explicitly.
\end{enumerate}
\end{theorem}

Theorem~\ref{trans} was used in \cite{giudici} to classify all
innately transitive groups having no fixed-point-free elements of
prime order, settling the polycirculant graph conjecture for such
groups. Another application of Theorems~\ref{factns} and
\ref{trans} is the following factorisation theorem that was used
in the proof of Theorem~\ref{qpinclns}.  It implies in particular
that, if $G$ is quasiprimitive of degree $n$ with nonabelian and
non-simple socle, then $S_n$ and possibly $A_n$ are the only
almost simple over-groups of $G$.

\begin{theorem}{\rm \cite[Theorem 1.4]{qpinclns2}}
Let $T, S$ be finite nonabelian simple groups such that $T$ has
proper subgroups $A, B$ with $T=AB$ and $A= S^\ell$ for some
$\ell\geq2$. Then $T=A_n$, $B=A_{n-1}$, where $n=|T:B|$, and $A$
is a transitive group of degree $n$.
\end{theorem}

Finally we note that Theorem~\ref{degs} is based on the following
result about indices of subgroups of finite simple groups.

\begin{theorem}{\rm \cite{CNT,PraegerShalev02}}
For a positive real number $x$, the proportion of integers $n\leq
x$ of the form $n=|T:M|$, where $T$ is a nonabelian simple group
and $M$ is either a maximal subgroup or a proper subgroup, and
$(T,M)\ne (A_n,A_{n-1})$, is at most $(1+o(1))c/\log x$, where
$c=1$ or $c=\sum_{d=1}^\infty\frac{1}{d\phi(d)}$ respectively.
\end{theorem}

We have presented a framework for studying finite permutation
groups by identifying and analysing their basic components.  The
impetus for extending the theory beyond primitive groups came from
the need for an appropriate theory of basic permutation groups for
combinatorial applications.  Developing this theory required the
answers to specific questions about simple groups, and the power
of the theory is largely due to its use of the finite simple group
classification.

\label{lastpage}


\begin{thebibliography}{aa}
\bibitem{asch}
M. Aschbacher, On the maximal subgroups of the finite classical
groups, \emph{Invent. Math.} {\bf 76} (1984), 469--514.
\bibitem{qpinclns}
R. Baddeley and C. E. Praeger, On primitive overgroups of
quasiprimitive permutation groups, Research Report No. 2002/3, U.
Western Australia, 2002.
\bibitem{BambergPraeger}
J. Bamberg and C. E. Praeger, Finite permutation groups with a
transitive minimal normal subgroup, preprint, 2002.
\bibitem{Cameron81}
P. J. Cameron, Finite permutation groups and finite simple groups,
\emph{Bull. London Math. Soc.} {\bf 13} (1981), 1--22.
\bibitem{CNT} P. J. Cameron, P. M. Neumann and D. N. Teague, On
the degrees of primitive permutation groups, {\it Math. Z.} {\bf
180} (1982), 141--149.
\bibitem{CameronPraegeretal}
P. J. Cameron, C. E. Praeger, J. Saxl, and G. M. Seitz, On the
Sims conjecture and distance transitive graphs, {\it Bull. London
Math. Soc.} {\bf 15} (1983), 499--506.
\bibitem{fhp}
X. G. Fang, G. Havas, and C. E. Praeger, On the automorphism
groups of quasiprimitive almost simple graphs, \emph{J. Algebra}
{\bf 222} (1999), 271--283.
\bibitem{cayley}
X. G. Fang, C. E. Praeger and J. Wang, On the automorphism groups
of Cayley graphs of finite simple groups, \emph{J. London Math.
Soc.} (to appear).
\bibitem{fgs}
M. D. Fried, R. Guralnick\ and\ J. Saxl, Schur covers and
Carlitz's conjecture, \emph{Israel J. Math.} {\bf 82} (1993),
157--225.
\bibitem{giudici}
M. Giudici, Quasiprimitive groups with no fixed point free
elements of prime order, \emph{J. London Math. Soc.}, (to appear).
\bibitem{giudici2}
M. Giudici, C. H. Li and C. E. Praeger, Analysing finite locally
$s$-arc-transitive graphs, in preparation.
\bibitem{sasha}
A. A. Ivanov, Distance-transitive graphs and their classification, in \emph{Investigations in algebraic theory of
combinatorial objects}, Kluwer, Dordrecht, 1994, 283--378.
\bibitem{ip}
A. A. Ivanov and C. E. Praeger, On finite affine $2$-arc
transitive graphs, \emph{European J. Combin.} {\bf 14} (1993),
421--444.
\bibitem{kovacs}
L. G. Kovacs, Primitive subgroups of wreath products in product
action, \emph{Proc. London Math. Soc.} (3) {\bf 58} (1989),
306--322.
\bibitem{li1}
C. H. Li, Finite s-arc transitive graphs of prime-power order,
\emph{Bull. London Math. Soc.} {\bf 33} (2001), 129-137.
\bibitem{li2}
C. H. Li, On finite $s$-arc transitive graphs of odd order,
\emph{J. Combin. Theory Ser. B} {\bf 81} (2001), 307-317.
\bibitem{li3}
C. H. Li, The finite vertex-primitive and vertex-biprimitive
$s$-transitive graphs for $s\geq 4$, \emph{Trans. Amer. Math.
Soc.} {\bf 353} (2001), 3511--3529.
\bibitem{maxan}
M. W. Liebeck, C. E. Praeger and J. Saxl, A classification of the
maximal subgroups of the finite alternating and symmetric groups,
\emph{Proc. London Math. Soc.} {\bf 55} (1987), 299--330.
\bibitem{factns}
M. W. Liebeck, C. E. Praeger and J. Saxl, The maximal
factorisations of the finite simple groups and their automorphism
groups, {\it Mem. Amer. Math. Soc.} No. 432, Vol. 86 (1990),
1--151.
\bibitem{max+}
M. W. Liebeck, C. E. Praeger and J. Saxl, On factorisations of
almost simple groups, {\it J. Algebra} {\bf 185} (1996), 409--419.
\bibitem{trans}
M. W. Liebeck, C. E. Praeger and J. Saxl, Transitive subgroups of
primitive permutation groups, {\it J. Algebra} {\bf 234} (2000),
291--361.
\bibitem{LiebeckPraegeretal02}
M. W. Liebeck, C. E. Praeger and J. Saxl, Primitive permutation
groups with a common suborbit, and edge-transitive graphs, {\it
Proc. London Math. Soc.} (3) {\bf 84} (2002), 405--438.
\bibitem{inclns}
C. E. Praeger, The inclusion problem for finite primitive
permutation groups, {\it Proc. London Math. Soc.} (3) {\bf 60}
(1990), 68--88.
\bibitem{2arcqp}
C. E. Praeger, \newblock An {O}'{N}an-{S}cott theorem for finite
quasiprimitive permutation groups and an application to 2-arc
transitive graphs, \newblock {\em J. London Math. Soc.} (2) 47
(1993), 227--239.
\bibitem{bcc}
C. E. Praeger, Quasiprimitive graphs. In \emph{Surveys in combinatorics, 1997 (London)}, 65--85, Cambridge
University Press, Cambridge, 1997.
\bibitem{qpinclns2}
C. E. Praeger, Quotients and inclusions of finite quasiprimitive
permutation groups, Research Report No. 2002/05, University of
Western Australia, 2002.
\bibitem{seminorm}
C. E. Praeger, Seminormal and subnormal subgroup lattices for
transitive permutation groups, in preparation.
\bibitem{psy}
C. E. Praeger, J. Saxl and K. Yokoyama, Distance transitive graphs
and finite simple groups, \emph{Proc. London Math. Soc.} (3) {\bf
55} (1987), 1--21.
\bibitem{PraegerShalev01}
C. E. Praeger and A. Shalev, Bounds on finite quasiprimitive
permutation groups, {\it J. Austral. Math. Soc.} {\bf 71} (2001),
243--258.
\bibitem{PraegerShalev02}
C. E. Praeger and A. Shalev, Indices of subgroups of finite simple
groups and quasiprimitive permutation groups, preprint, 2002.
\bibitem{vbc}
J.  van Bon and A. M. Cohen, Prospective classification of distance-trans\-it\-ive graphs, in \emph{Combinatorics
'88 (Ravello)}, Mediterranean, Rende, 1991, 25--38.
\end{thebibliography}
\end{document}